# Combined High-Dimensional and Reduced-Order Modeling Based on an Iterative Domain Decomposition Method

Taiji Saito *, Tomoshi Miyamura ** and Yasunori Yusa ***

* Department of Computer Science, Graduate School of Engineering, Nihon University

1 Nakagawara, Tokusada, Tamura-machi, Koriyama, Fukushima 963-8642, Japan

** Department of Computer Science, College of Engineering, Nihon University

1 Nakagawara, Tokusada, Tamura-machi, Koriyama, Fukushima 963-8642, Japan

E-mail: miyamura.tomoshi@nihon-u.ac.jp

*** Department of Mechanical and Intelligent Systems Engineering, Graduate School of Informatics and Engineering, The University of Electro-Communications

1-5-1 Chofugaoka, Chofu, Tokyo 182-8585, Japan

**Abstract**

A method is developed for structural analysis in which parts of an analysis domain are modeled by reduced-order models (ROMs) and the remaining part is modeled by a high-dimensional model (HDM). The multiple ROMs and HDM are combined based on an iterative domain decomposition method (DDM) with an arbitrary number of non-overlapping subdomains. A diagonal preconditioner is derived for the iterative DDM. The method is implemented based on a framework of the existing DDM-based analysis code. The ROM is implemented on GPU. Snapshot generation methods for the ROM are also proposed. Static structural analyses of nuclear power plant models are conducted using the combined multiple ROMs and HDM with multiple-subdomains model. The analyses are performed in parallel by the hybrid MPI-CUDA implementation. The convergence property of the iterative DDM is improved by the preconditioner.

***Keywords*** : Reduced-Order Model, ROM, High-Dimensional Model, HDM, Domain Decomposition Method, Iterative Solver

## 1. Introduction

When the structural analysis of an artificial structure with a complex geometry is performed using the finite element method (FEM), a large number of elements is required to represent the geometry accurately. Furthermore, finer mesh discretization is required to improve the analysis accuracy. Owing to advances in parallel computing techniques and computer performance, high-fidelity analyses using such large-scale high-dimensional models (HDMs) have become feasible (Yoshimura et al., 2019; Miyamura et al., 2019). However, since a large number of analyses must be performed for engineering studies, computational acceleration is required.

In recent years, projection-based reduced-order models (ROMs) have been proposed as a means of accelerating computations by reducing the dimension of the analysis model while minimizing the loss of analysis accuracy (Farhat et al., 2014; Tezaur et al., 2022). These methods have been applied to various problems in computational mechanics. In these methods, solutions under various parameters are first computed using an HDM and are collected as snapshots. Next, useful reduced-order basis (ROB) is extracted by applying the proper orthogonal decomposition (POD) or the singular value decomposition (SVD) to the matrix composed of the snapshots. Finally, a ROM can be constructed by applying the Petrov-Galerkin projection to the governing equations in fluid flow analysis or the Galerkin projection to the governing equations in structural analysis using the extracted ROB.

When the dimension of an HDM is large, the dimension of the snapshots also becomes large. Therefore, the Local-POD method, which constructs local ROBs based on domain decomposition, was proposed by Baiges et al. (2013). In this method, the analysis domain is divided into overlapping subdomains, and a ROM is constructed for each subdomain. These ROMs are then used to construct a ROM for the entire problem. Furthermore, an FOM-ROM domain decomposition method, in which part of the analysis domain is modeled by a ROM and the remaining part is modeled by a full-order model (FOM), i.e., an HDM, was also proposed. These methods have been applied to fluid flow analysis based on the Navier-Stokes equations. For the Local-POD method, Shintate et al. (2025) proposed a parallel solution method for linear systems based on a coarser domain decomposition than that used for the local-POD method. In this method, the coarser domain decomposition is obtained by partitioning a graph that represents the relationships among the local PODs obtained by the original domain decomposition.

Similar methods have also been proposed for structural analysis (Kerfriden et al., 2011; Radermacher et al., 2014; Corigliano et al., 2015). In these methods, highly nonlinear region in nonlinear structural analysis is modeled by HDM, whereas nearly elastic region is modeled by ROM, and these models are coupled based on the substructuring method. The analysis domain is divided into non-overlapping two subdomains, that is, a ROM subdomain and a HDM subdomain. In Kerfriden et al. (2011), quasi-brittle fracture region was modeled by HDM. Adaptive methods for identifying an elastic-plastic region to be modeled by HDM were proposed by Radermacher et al. (2014) and Corigliano et al. (2015). Corigliano et al. (2015) introduced interface stiffnesses at the boundaries between subdomains to couple the subdomains. Kerfriden et al. (2011) and Radermacher et al. (2014) coupled the subdomains by sharing the degrees of freedom (DOFs) on the interface between the ROM and HDM subdomains, and their methods were based on the so-called primal substructuring method. Kerfriden et al. (2011) also proposed a method in which unknowns in the ROM subdomain are solved using a direct solver and the remaining unknowns are solved by the preconditioned conjugate gradient method.

In the present study, a method is proposed for structural analysis in which parts of an analysis domain are modeled by ROMs and the remaining part is modeled by an HDM. In the previous studies by Kerfriden et al. (2011), Radermacher et al. (2014), and Corigliano et al. (2015), the analysis domain is divided into two subdomains, i.e., a subdomain modeled by a nonlinear HDM and a subdomain modeled by a ROM. In the present study, multiple ROMs and HDM are combined based on an iterative domain decomposition method (DDM) using non-overlapping subdomains. The HDM part is divided into an arbitrary number of subdomains. The linear algebraic equations derived from the implicit finite element method are solved by an iterative solver such as the conjugate gradient (CG) method, similar to the approach proposed by Kerfriden et al. (2011). The subdomains to be modeled by ROMs are first defined, and the remaining HDM domain is then divided into non-overlapping subdomains. The proposed method can be incorporated into the framework of existing parallel iterative solvers and iterative primal substructuring methods based on non-overlapping domain decomposition (Yoshimura et al., 2002; ADVENTURE Project, 2026; Miyamura et al., 2023).

Furthermore, in each iteration of the CG method, the computation for a subdomain modeled by ROM becomes a local problem with Dirichlet boundary conditions imposed on the subdomain boundaries. Therefore, snapshots with parameters that represent the Dirichlet boundary conditions are required for the ROMs, and methods for generating such snapshots are also proposed. Such systematic methods for snapshot generation have not been reported in previous studies. In addition, a preconditioning method for the CG method used in the present method is also proposed. Many artificial structures are composed of a large number of components. A detailed finite element analysis of such assembled structures is one of the problems to which the proposed method can be effectively applied. Components of engineering interest or components with strong nonlinearities are modeled by HDM, whereas the remaining components are modeled by ROMs.

The remainder of this paper is organized as follows. In Section 2, an overview of the ROM based on the Galerkin projection for structural analysis is presented. In Section 3, the formulation of the domain decomposition method is briefly described. In Section 4, the solution procedure based on the preconditioned conjugate gradient method is described. In Section 5, the method for combining ROMs and HDM based on an iterative DDM (HDM–ROM method) is proposed. In Section 6, a method for generating snapshots for the ROMs used in the HDM–ROM method is investigated. Section 7 describes the implementation of the proposed method, including the GPU implementation of the ROM. Numerical examples are presented in Section 8. Finally, conclusions are given in Section 9.

## 2. Overview of Galerkin Projection-Based ROM for Structural Analysis

A ROM based on Galerkin projection is briefly described with reference to Farhat et al. (2014) and Tezaur et al. (2022). Static linear structural analysis based on the finite element method is considered. The equilibrium equation representing an HDM is expressed as follows:

$$\boldsymbol{K}(\boldsymbol{\mu})\boldsymbol{u}(\boldsymbol{\mu}) = \boldsymbol{f}(\boldsymbol{\mu}) \tag{1}$$

where $\boldsymbol{\mu}$ denotes a parameter vector, $\boldsymbol{K}(\boldsymbol{\mu})$, $\boldsymbol{u}(\boldsymbol{\mu})$, and $\boldsymbol{f}(\boldsymbol{\mu})$ denote, the stiffness matrix, the nodal displacement vector, and the load vector. In structural analysis, the parameters $\boldsymbol{\mu}$ may include material properties, the geometry of the analysis domain, and boundary conditions. Eq. (1) is rewritten as follows:

$$\boldsymbol{r}(\boldsymbol{u}(\boldsymbol{\mu});\boldsymbol{\mu}) \equiv \boldsymbol{f}(\boldsymbol{\mu}) - \boldsymbol{K}(\boldsymbol{\mu})\boldsymbol{u}(\boldsymbol{\mu}) = \boldsymbol{0} \tag{2}$$

The solution $\boldsymbol{u}(\boldsymbol{\mu})$ is approximated by an affine subspace:

$$\boldsymbol{u}(\boldsymbol{\mu}) \approx \boldsymbol{u}_{\mathrm{ref}} + \boldsymbol{V}\boldsymbol{y}(\boldsymbol{\mu}) \tag{3}$$

where $\boldsymbol{V}$ denotes the reduced-order basis (ROB), and $\boldsymbol{u}_{ref}$ is prescribed appropriately. By using these quantities, the governing equation of the ROM based on the Galerkin projection of Eq. (2) is obtained as follows:

$$\hat{\boldsymbol{r}}(\boldsymbol{y}(\boldsymbol{\mu});\boldsymbol{\mu}) \equiv \boldsymbol{V}^{\mathrm{T}}\boldsymbol{r}(\boldsymbol{u}_{ref} + \boldsymbol{V}\boldsymbol{y}(\boldsymbol{\mu});\boldsymbol{\mu}) = \boldsymbol{0} \tag{4}$$

The procedure for obtaining the ROB $\boldsymbol{V}$ is described below. First, analyses are performed using the HDM for $m$ sets of $\boldsymbol{\mu}$ (denoted by $\boldsymbol{\mu}_i$, $i = 1, \dots, m$). The corresponding solution vector $\boldsymbol{u}(\boldsymbol{\mu}_i)$ is referred to as a snapshot. By arranging the snapshots $\boldsymbol{u}(\boldsymbol{\mu}_i)$, the following matrix is defined:

$$\boldsymbol{S} = \left[\boldsymbol{u}(\boldsymbol{\mu}_1) - \boldsymbol{u}_{ref}, \dots, \boldsymbol{u}(\boldsymbol{\mu}_k) - \boldsymbol{u}_{ref}\right] \tag{5}$$

$\boldsymbol{u}_{ref}$ in Eqs. (3) and (5) is set to the average of the snapshots $\boldsymbol{u}(\boldsymbol{\mu}_i)$ ($i = 1, \dots, m$) or another appropriate value. In the present study, it is set to $\boldsymbol{0}$. Next, SVD of $\boldsymbol{S}^{(k)}$ is obtained as follows:

$$\boldsymbol{S} = \boldsymbol{U}_{\mathrm{SVD}}\boldsymbol{\Sigma}_{\mathrm{SVD}}\boldsymbol{W}_{\mathrm{SVD}} \tag{6}$$

$\boldsymbol{V}$ is constructed from the first $r$ columns of $\boldsymbol{U}_{\mathrm{SVD}}$, which are used as the ROB vectors. The number of ROB vectors $r$ is determined as the smallest integer satisfying the following inequality:

$$E(r) = \sum_{i=1}^{r}\sigma_i \Big/ \sum_{i=1}^{M}\sigma_i \geq \eta \tag{7}$$

where $\sigma_i$ denotes the $i$-th singular value, $M$ is the number of snapshots, and $\eta$ is the prescribed threshold for the cumulative energy contribution ratio $E(r)$. By substituting Eq. (3) into Eq. (2), and subsequently substituting the resulting equation into Eq. (4), the following equation is obtained:

$$\hat{\boldsymbol{r}}(\boldsymbol{y}(\boldsymbol{\mu});\boldsymbol{\mu}) \equiv \boldsymbol{V}^{\mathrm{T}}\left(\boldsymbol{f}(\boldsymbol{\mu}) - \boldsymbol{K}(\boldsymbol{\mu})\boldsymbol{V}\boldsymbol{y}(\boldsymbol{\mu})\right) = \boldsymbol{0} \tag{8}$$

By rearranging Eq. (8), the governing equation of the Galerkin projection-based ROM (hereinafter referred to as the ROM) is obtained as follows:

$$\left(\boldsymbol{V}^{\mathrm{T}}\boldsymbol{K}(\boldsymbol{\mu})\boldsymbol{V}\right)\boldsymbol{y}(\boldsymbol{\mu})=\boldsymbol{V}^{\mathrm{T}}\boldsymbol{f}(\boldsymbol{\mu}) \tag{9}$$

By solving Eq. (9) for $\boldsymbol{y}(\boldsymbol{\mu})$ and substituting the result into Eq. (3), an approximation of the nodal displacement vector can be computed. Since $\boldsymbol{y}(\boldsymbol{\mu})$ is an unknown variable determined by solving Eq. (9), it is hereinafter denoted by $\boldsymbol{y}$.

## 3. Domain Decomposition Method

### 3.1 Non-Overlapping Domain Decomposition Method

The equilibrium equation given in Eq. (1) is solved by a domain decomposition method with an iterative solver (iterative DDM). First, the analysis domain is divided into $N$ non-overlapping subdomains. For subdomain $k$, the stiffness matrix, nodal displacement vector, and nodal load vector are denoted by $\boldsymbol{K}^{(k)}$, $\boldsymbol{u}^{(k)}$, and $\boldsymbol{f}^{(k)}$, respectively. The stiffness matrix $\boldsymbol{K}$ is obtained by assembling the stiffness matrices $\boldsymbol{K}^{(k)}$ of the subdomains as follows:

$$\boldsymbol{K}=\sum_{k=1}^{N}\boldsymbol{N}^{(k)}\boldsymbol{K}^{(k)}\boldsymbol{N}^{(k)\mathrm{T}} \tag{10}$$

where $\boldsymbol{N}^{(k)}$ denotes a Boolean matrix that maps the local DOFs of a subdomain to the global DOFs. The sharing of nodes by multiple subdomains at the interfaces between subdomains is also represented by $\boldsymbol{N}^{(k)}$. The equilibrium equation represented by Eq. (1) can then be written as follows:

$$\sum_{k=1}^{N}\boldsymbol{N}^{(k)}\boldsymbol{K}^{(k)}\boldsymbol{u}^{(k)}=\sum_{k=1}^{N}\boldsymbol{N}^{(k)}\boldsymbol{f}^{(k)} \tag{11}$$

Note that $\boldsymbol{u}^{(k)}=\boldsymbol{N}^{(k)\mathrm{T}}\boldsymbol{u}$. Note also that when an iterative method such as the CG method is parallelized based on the domain-decomposed equilibrium equation given in Eq. (11), such a solver is referred to as a parallel iterative solver.

### 3.2 Substructuring-Based Domain Decomposition Method

For subdomain $k$, $\boldsymbol{K}^{(k)}$, $\boldsymbol{u}^{(k)}$, and $\boldsymbol{f}^{(k)}$ are partitioned into the DOFs in the interior of the subdomain and those on the interfaces between subdomains, as follows:

$$\boldsymbol{K}^{(k)}=\begin{bmatrix}\boldsymbol{K}_{\mathrm{II}}^{(k)} & \boldsymbol{K}_{\mathrm{IB}}^{(k)}\\ \boldsymbol{K}_{\mathrm{BI}}^{(k)} & \boldsymbol{K}_{\mathrm{BB}}^{(k)}\end{bmatrix},\quad \boldsymbol{u}^{(k)}=\begin{bmatrix}\boldsymbol{u}_{\mathrm{I}}^{(k)}\\ \boldsymbol{u}_{\mathrm{B}}^{(k)}\end{bmatrix},\quad \boldsymbol{f}^{(k)}=\begin{bmatrix}\boldsymbol{f}_{\mathrm{I}}^{(k)}\\ \boldsymbol{f}_{\mathrm{B}}^{(k)}\end{bmatrix} \tag{12}$$

where the subscript I denotes the DOFs in the interior of a subdomain, and the subscript B denotes the DOFs on the interfaces between subdomains. The total number of DOFs on all interfaces between subdomains is denoted by $n_{\mathrm{B}}$. In addition, the number of interface DOFs associated with subdomain $k$ is denoted by $n_{\mathrm{B}}^{(k)}$.

The equilibrium associated with the interior DOFs of the subdomains is satisfied independently in each subdomain as follows:

$$\boldsymbol{K}_{\mathrm{II}}^{(k)}\boldsymbol{u}_{\mathrm{I}}^{(k)}+\boldsymbol{K}_{\mathrm{IB}}^{(k)}\boldsymbol{u}_{\mathrm{B}}^{(k)}=\boldsymbol{f}_{\mathrm{I}}^{(k)} \tag{13}$$

On the other hand, the equilibrium of the interface DOFs between subdomains is expressed as follows.

$$\sum_{k=1}^{N}\boldsymbol{N}_{\mathrm{B}}^{(k)}\boldsymbol{K}_{\mathrm{BI}}^{(k)}\boldsymbol{u}_{\mathrm{I}}^{(k)}+\sum_{k=1}^{N}\boldsymbol{N}_{\mathrm{B}}^{(k)}\boldsymbol{K}_{\mathrm{BB}}^{(k)}\boldsymbol{u}_{\mathrm{B}}^{(k)}=\sum_{k=1}^{N}\boldsymbol{N}_{\mathrm{B}}^{(k)}\boldsymbol{f}_{\mathrm{B}}^{(k)} \tag{14}$$

where $\boldsymbol{N}_{\mathrm{B}}^{(k)}$ denotes a Boolean matrix that maps the interface DOFs associated with subdomain $k$ to the interface

DOFs of all subdomains. Equation (13) can be solved for the interior nodal displacements as follows:

$$\boldsymbol{u}_{\mathrm{I}}^{(k)} = \left(\boldsymbol{K}_{\mathrm{II}}^{(k)}\right)^{-1} \left( \boldsymbol{f}_{\mathrm{I}}^{(k)} - \boldsymbol{K}_{\mathrm{IB}}^{(k)} \boldsymbol{u}_{\mathrm{B}}^{(k)} \right) \tag{15}$$

By substituting Eq. (15) into Eq. (14) and moving the known terms to the right-hand side, the following equation is obtained:

$$\sum_{k=1}^{N} \boldsymbol{N}_{\mathrm{B}}^{(k)} \left\{ \boldsymbol{K}_{\mathrm{BB}}^{(k)} - \boldsymbol{K}_{\mathrm{BI}}^{(k)} \left(\boldsymbol{K}_{\mathrm{II}}^{(k)}\right)^{-1} \boldsymbol{K}_{\mathrm{IB}}^{(k)} \right\} \boldsymbol{u}_{\mathrm{B}}^{(k)} = \sum_{k=1}^{N} \boldsymbol{N}_{\mathrm{B}}^{(k)} \left\{ \boldsymbol{f}_{\mathrm{B}}^{(k)} - \boldsymbol{K}_{\mathrm{BI}}^{(k)} \left(\boldsymbol{K}_{\mathrm{II}}^{(k)}\right)^{-1} \boldsymbol{f}_{\mathrm{I}}^{(k)} \right\} \tag{16}$$

The local Schur complement $\boldsymbol{S}^{(k)}$ for subdomain $k$ is defined as follows:

$$\boldsymbol{S}^{(k)} = \boldsymbol{K}_{\mathrm{BB}}^{(k)} - \boldsymbol{K}_{\mathrm{BI}}^{(k)} \left(\boldsymbol{K}_{\mathrm{II}}^{(k)}\right)^{-1} \boldsymbol{K}_{\mathrm{IB}}^{(k)} \tag{17}$$

The corresponding right-hand-side vector is also defined as follows:

$$\tilde{\boldsymbol{f}}_{\mathrm{B}}^{(k)} = \boldsymbol{f}_{\mathrm{B}}^{(k)} - \boldsymbol{K}_{\mathrm{BI}}^{(k)} \left(\boldsymbol{K}_{\mathrm{II}}^{(k)}\right)^{-1} \boldsymbol{f}_{\mathrm{I}}^{(k)} \tag{18}$$

By using Eqs. (17) and (18), Eq. (16) can be expressed as follows:

$$\sum_{k=1}^{N} \boldsymbol{N}_{\mathrm{B}}^{(k)} \boldsymbol{S}^{(k)} \boldsymbol{u}_{\mathrm{B}}^{(k)} = \sum_{k=1}^{N} \boldsymbol{N}_{\mathrm{B}}^{(k)} \tilde{\boldsymbol{f}}_{\mathrm{B}}^{(k)} \tag{19}$$

By assembling the local Schur complements $\boldsymbol{S}^{(k)}$ for all subdomains, the global Schur complement $\boldsymbol{S}$ is obtained as follows:

$$\boldsymbol{S} = \sum_{k=1}^{N} \boldsymbol{N}_{\mathrm{B}}^{(k)} \boldsymbol{S}^{(k)} \boldsymbol{N}_{\mathrm{B}}^{(k)\mathrm{T}} \tag{20}$$

### 3.3 Combination of Non-Substructuring and Substructuring-Based DDMs

A method for combining the non-substructuring and substructuring-based DDMs is described. The number of subdomains without condensation of the interior DOFs is $N_{\mathrm{K}}$, and the number of substructuring-based subdomains is $N_{\mathrm{S}}$. In this case, the total number of DOFs is reduced from that of the original problem by the number of interior DOFs condensed in the substructuring-based subdomains. The reduced DOFs are called the modified global DOFs.. Accordingly, the numbers of rows of the Boolean matrices $\boldsymbol{N}^{(k)}$ and $\boldsymbol{N}_{\mathrm{B}}^{(k)}$ are modified so that they correspond to the modified global DOFs. Let $\tilde{\boldsymbol{u}}$ denote the nodal displacement vector associated with the modified global DOFs. Then, $\boldsymbol{u}_{\mathrm{B}}^{(k)} = \boldsymbol{N}_{\mathrm{B}}^{(k)\mathrm{T}} \tilde{\boldsymbol{u}}$ and $\boldsymbol{u}^{(k)} = \boldsymbol{N}^{(k)\mathrm{T}} \tilde{\boldsymbol{u}}$ hold. The equilibrium equation is expressed as follows:

$$\sum_{k=1}^{N_{\mathrm{K}}} \boldsymbol{N}^{(k)} \boldsymbol{K}^{(k)} \boldsymbol{u}^{(k)} + \sum_{k=1}^{N_{\mathrm{B}}} \boldsymbol{N}_{\mathrm{B}}^{(k)} \boldsymbol{S}^{(k)} \boldsymbol{u}_{\mathrm{B}}^{(k)} = \sum_{k=1}^{N_{\mathrm{K}}} \boldsymbol{N}^{(k)} \boldsymbol{f}^{(k)} + \sum_{k=1}^{N_{\mathrm{B}}} \boldsymbol{N}_{\mathrm{B}}^{(k)} \tilde{\boldsymbol{f}}_{\mathrm{B}}^{(k)} \tag{21}$$

Equation (21) can be rewritten as follows by using $\tilde{\boldsymbol{u}}$:

$$\left( \sum_{k=1}^{N_{\mathrm{K}}} \boldsymbol{N}^{(k)} \boldsymbol{K}^{(k)} \boldsymbol{N}^{(k)\mathrm{T}} + \sum_{k=1}^{N_{\mathrm{B}}} \boldsymbol{N}_{\mathrm{B}}^{(k)} \boldsymbol{S}^{(k)} \boldsymbol{N}_{\mathrm{B}}^{(k)\mathrm{T}} \right) \tilde{\boldsymbol{u}} = \sum_{k=1}^{N_{\mathrm{K}}} \boldsymbol{N}^{(k)} \boldsymbol{f}^{(k)} + \sum_{k=1}^{N_{\mathrm{B}}} \boldsymbol{N}_{\mathrm{B}}^{(k)} \tilde{\boldsymbol{f}}_{\mathrm{B}}^{(k)} \tag{22}$$

By defining $\tilde{\boldsymbol{S}} = \sum_{k=1}^{N_{\mathrm{K}}} \boldsymbol{N}^{(k)} \boldsymbol{K}^{(k)} \boldsymbol{N}^{(k)\mathrm{T}} + \sum_{k=1}^{N_{\mathrm{B}}} \boldsymbol{N}_{\mathrm{B}}^{(k)} \boldsymbol{S}^{(k)} \boldsymbol{N}_{\mathrm{B}}^{(k)\mathrm{T}}$ and $\tilde{\boldsymbol{f}} = \sum_{k=1}^{N_{\mathrm{K}}} \boldsymbol{N}^{(k)} \boldsymbol{f}^{(k)} + \sum_{k=1}^{N_{\mathrm{B}}} \boldsymbol{N}_{\mathrm{B}}^{(k)} \tilde{\boldsymbol{f}}_{\mathrm{B}}^{(k)}$, Eq. (22) can be written as follows:

$$\tilde{\boldsymbol{S}} \boldsymbol{u} = \tilde{\boldsymbol{f}} \tag{23}$$

## 4. Solution Procedure of the CG Method within the Framework of the DDM

Equation (23) is solved using the preconditioned CG method. Algorithm 1 provides the procedure for the CG method, where $\boldsymbol{M}^{-1}$ denotes a preconditioning matrix. During the procedure of a CG iteration for Eq. (23), the product of $\tilde{\boldsymbol{S}}$ and the search direction vector $\boldsymbol{p}^i$ is computed. This computation includes the evaluation of the product of the local Schur complement $\boldsymbol{S}^{(k)}$ and the vector $\boldsymbol{p}_{\mathrm{B}}^{(k)}$, as shown below, in order to obtain the nodal reaction force vector $\boldsymbol{r}_{\mathrm{B}}^{(k)}$. Here, $\boldsymbol{p}_{\mathrm{B}}^{(k)}$ denotes the block of $\boldsymbol{p}^i$ corresponding to the interface DOFs of subdomain $k$.

$$\boldsymbol{r}_{\mathrm{B}}^{(k)} = \boldsymbol{S}^{(k)} \boldsymbol{p}_{\mathrm{B}}^{(k)} \tag{24}$$

The reaction forces are also computed for the subdomains in which the interior DOFs are not condensed, and the reaction forces represented by Eq. (24) are assembled together with those of all other subdomains.

In practice, it is more efficient to solve a problem in which the interface DOFs between subdomains are prescribed as Dirichlet boundary conditions in addition to the original Dirichlet boundary conditions for prescribed displacements (referred to as a local Dirichlet problem), as shown below, rather than explicitly constructing $\boldsymbol{S}^{(k)}$ and performing the computation using Eq. (24).

$$\begin{bmatrix} \boldsymbol{K}_{\mathrm{II}}^{(k)} & \boldsymbol{K}_{\mathrm{IB}}^{(k)} \\ \boldsymbol{K}_{\mathrm{BI}}^{(k)} & \boldsymbol{K}_{\mathrm{BB}}^{(k)} \end{bmatrix} \begin{bmatrix} \boldsymbol{u}_{\mathrm{I}}^{(k)} \\ \boldsymbol{p}_{\mathrm{B}}^{(k)} \end{bmatrix} = \begin{bmatrix} \boldsymbol{f}_{\mathrm{I}}^{(k)} \\ \boldsymbol{r}_{\mathrm{B}}^{(k)} \end{bmatrix} \tag{25}$$

From Eq. (25), the following linear system with $\boldsymbol{K}_{\mathrm{II}}^{(k)}$ as the coefficient matrix is obtained.

$$\boldsymbol{K}_{\mathrm{II}}^{(k)} \boldsymbol{u}_{\mathrm{I}}^{(k)} = -\boldsymbol{K}_{\mathrm{IB}}^{(k)} \boldsymbol{p}_{\mathrm{B}}^{(k)} + \boldsymbol{f}_{\mathrm{I}}^{(k)} \tag{26}$$

After solving this system using a linear solver for $\boldsymbol{u}_{\mathrm{I}}^{(k)}$ in the same manner as Eq. (15), the reaction force $\boldsymbol{r}_{\mathrm{B}}^{(k)}$ is obtained as follows.

$$\boldsymbol{r}_{\mathrm{B}}^{(k)} = \boldsymbol{K}_{\mathrm{BI}}^{(k)} \boldsymbol{u}_{\mathrm{I}}^{(k)} + \boldsymbol{K}_{\mathrm{BB}}^{(k)} \boldsymbol{p}_{\mathrm{B}}^{(k)} \tag{27}$$

| | |
|---|---|
| Initialization : $\begin{aligned} &\boldsymbol{u}^0 \text{ : arbitrary initial vector} \\ &\boldsymbol{g}^0 = \boldsymbol{p} - \tilde{\boldsymbol{S}}\boldsymbol{u}^0 \\ &\boldsymbol{p}^0 = \boldsymbol{M}^{-1}\boldsymbol{g}^0 \end{aligned}$ | Iteration : $\begin{aligned} &\alpha = \frac{\left(\boldsymbol{M}^{-1}\boldsymbol{g}^i\right)^{\mathrm{T}}\boldsymbol{g}^i}{\boldsymbol{p}^{i\mathrm{T}}\left(\tilde{\boldsymbol{S}}\boldsymbol{p}^i\right)} \\ &\boldsymbol{u}^{i+1} = \boldsymbol{u}^i + \alpha\boldsymbol{p}^i \\ &\boldsymbol{g}^{i+1} = \boldsymbol{g}^i - \alpha\tilde{\boldsymbol{S}}\boldsymbol{p}^i \\ &\beta = \frac{\left(\boldsymbol{M}^{-1}\boldsymbol{g}^{i+1}\right)^{\mathrm{T}}\boldsymbol{g}^{i+1}}{\left(\boldsymbol{M}^{-1}\boldsymbol{g}^i\right)^{\mathrm{T}}\boldsymbol{g}^i} \\ &\boldsymbol{p}^{i+1} = \boldsymbol{M}^{-1}\boldsymbol{g}^{i+1} + \beta\boldsymbol{p}^i \\ &\text{if } \left\|\boldsymbol{g}^{i+1}\right\| / \left\|\boldsymbol{g}^1\right\| < \varepsilon_{converge} \Rightarrow \text{stop iteration} \end{aligned}$ |

Algorithm 1 The procedure for the CG method.

## 5. Combined High-Dimensional and Reduced-Order Modeling Based on Iterative DDM

Figure 1 shows a simplified model of a nuclear power plant used as an illustrative example in Section 7. The structure is a typical assembly structure composed of components such as a building and small modular reactors (SMRs). These components play different roles in the structural system. When the structural behavior of the building is of interest, a detailed analysis, including the evaluation of damage and elastic-plastic behavior, is required only for the building. For the other components, it is sufficient that their coupled effects on the building be evaluated appropriately. Conversely when the structural behavior of the SMRs is of interest, the detailed behavior of the building is not of primary concern. The region to be analyzed in detail is modeled by an HDM, whereas the remaining regions are modeled by ROMs. By coupling the HDM and the multiple ROMs, high accuracy is achieved in the region of interest, while the remaining regions are analyzed efficiently. The HDM and the multiple ROMs are combined using the framework of the iterative DDM; that is, the ROMs are incorporated into the solution procedure of the CG method described in Section 4. The method and the corresponding analysis model are referred to as the **HDM–ROM method** and the **HDM–ROM model**, respectively, in the present paper. Figure 2 illustrates the HDM–ROM model when the structural behavior of the building is of interest.

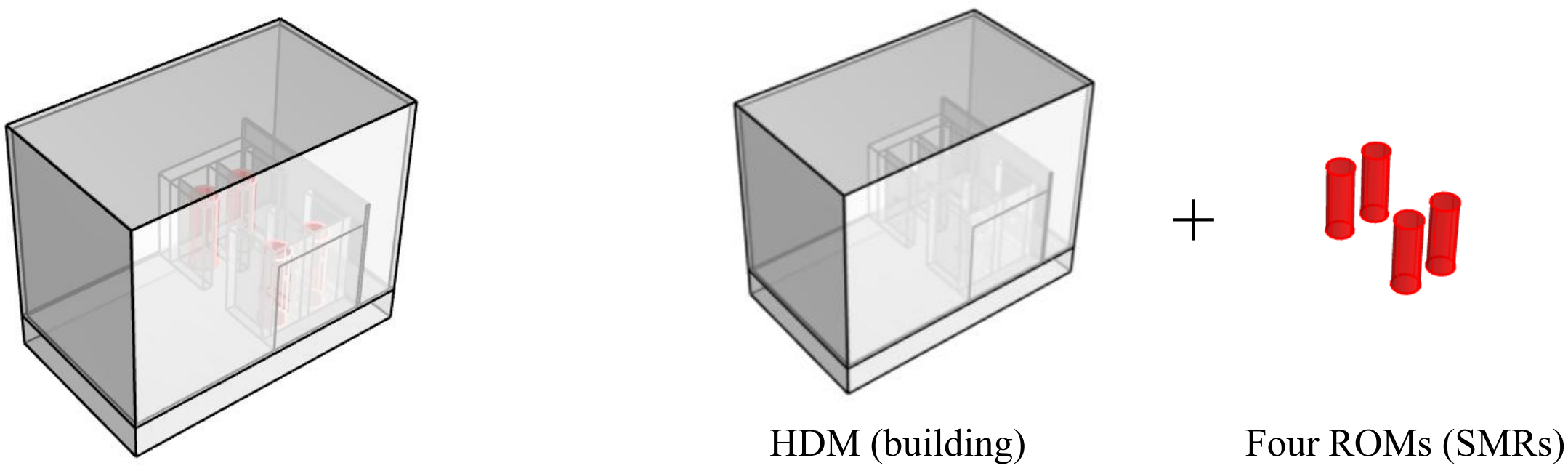

Fig. 1 A simplified model of a nuclear power plant with SMRs.

Fig. 2 The HDM–ROM model when the behavior of the building is of interest.

As shown in Eq. (26), the problem for a subdomain in the substructuring-based domain decomposition method is formulated as a Dirichlet boundary value problem for the subdomain. This problem is replaced by a ROM. Since the vector $\boldsymbol{p}_{\mathrm{B}}^{(k)}$ varies during the CG iterations, it is regarded as the parameter vector $\boldsymbol{\mu}$ of the ROM. In the present study, it is assumed that the load vector $\boldsymbol{f}_{\mathrm{I}}^{(k)}$ associated with the Neumann boundary conditions is $\boldsymbol{0}$ in the ROM region. Thus, Eq. (26) is rewritten as follows:

$$\boldsymbol{K}_{\mathrm{II}}^{(k)}\boldsymbol{u}_{\mathrm{I}}^{(k)} = -\boldsymbol{K}_{\mathrm{IB}}^{(k)}\boldsymbol{\mu} + \boldsymbol{f}_{\mathrm{I}}^{(k)} \tag{28}$$

Note that the number of parameters is $n_B^{(k)}$. Finally, the ROM equation corresponding to Eq. (9) is obtained as follows:

$$\boldsymbol{V}^T \boldsymbol{K}_{\mathrm{II}}^{(k)}(\boldsymbol{\mu}) \boldsymbol{V} \ \boldsymbol{y}^{(k)} = \boldsymbol{V}^T \left( -\boldsymbol{K}_{\mathrm{IB}}^{(k)} \boldsymbol{\mu} + \boldsymbol{f}_{\mathrm{I}}^{(k)} \right) \tag{29}$$

Equation (29) is solved for $\boldsymbol{y}^{(k)}$ using a linear solver. Then, $\boldsymbol{u}_{\mathrm{I}}^{(k)}$ is approximated as follows:

$$\boldsymbol{u}_{\mathrm{I}}^{(k)} \approx \boldsymbol{V}\boldsymbol{y}^{(k)} \tag{30}$$

Finally, the displacement vector $\boldsymbol{u}_{\mathrm{I}}^{(k)}$ obtained above is substituted into Eq. (27), and the reaction force $\boldsymbol{r}_{\mathrm{B}}^{(k)}$ is approximated. The above procedure is repeated within the CG iteration loop. Note that when the number of interior DOFs is sufficiently larger than that of the interface DOFs between subdomains, the computational cost of Eq. (27) is regarded as small.

The ROM approximation of the local Schur complement $\boldsymbol{S}^{(k)}$ defined by Eq. (17) can be described using $\left( \boldsymbol{V}^T \boldsymbol{K}_{\mathrm{II}}^{(k)}(\boldsymbol{\mu}) \boldsymbol{V} \right)^{-1}$. By denoting the ROM approximation of $\boldsymbol{S}^{(k)}$ by $\hat{\boldsymbol{S}}^{(k)}$, the following equation is obtained:

$$\boldsymbol{S}^{(k)} \approx \hat{\boldsymbol{S}}^{(k)} = \boldsymbol{K}_{\mathrm{BB}}^{(k)} - \boldsymbol{K}_{\mathrm{BI}}^{(k)} \boldsymbol{V} \left( \boldsymbol{V}^T \boldsymbol{K}_{\mathrm{II}}^{(k)}(\boldsymbol{\mu}) \boldsymbol{V} \right)^{-1} \boldsymbol{V}^T \boldsymbol{K}_{\mathrm{IB}}^{(k)} \tag{31}$$

Note that, as in the case of $\boldsymbol{S}^{(k)}$, it is not necessary to compute $\hat{\boldsymbol{S}}^{(k)}$ explicitly. Here, it is assumed that all subdomains whose interior DOFs are eliminated in Eq. (22) are solved using ROMs. The equilibrium equation is obtained by replacing $\boldsymbol{S}^{(k)}$ in Eq. (22) with $\hat{\boldsymbol{S}}^{(k)}$. The coefficient matrix on the left-hand side of Eq. (23) is redefined as $\tilde{\boldsymbol{S}} = \sum_{k=1}^{N_{\mathrm{K}}} \boldsymbol{N}^{(k)} \boldsymbol{K}^{(k)} \boldsymbol{N}^{(k)\mathrm{T}} + \sum_{k=1}^{N_{\mathrm{B}}} \boldsymbol{N}_{\mathrm{B}}^{(k)} \hat{\boldsymbol{S}}^{(k)} \boldsymbol{N}_{\mathrm{B}}^{(k)\mathrm{T}}$.

The diagonal entries of the coefficient matrix for diagonal scaling preconditioning are given as follows:

$$\boldsymbol{M}_{\mathrm{Diag}} = \sum_{k=1}^{N_{\mathrm{K}}} \boldsymbol{N}^{(k)} \mathrm{diag}\left( \boldsymbol{K}^{(k)} \right) \boldsymbol{N}^{(k)\mathrm{T}} + \sum_{k=1}^{N_{\mathrm{B}}} \boldsymbol{N}_{\mathrm{B}}^{(k)} \mathrm{diag}\left( \boldsymbol{S}^{(k)} \right) \boldsymbol{N}_{\mathrm{B}}^{(k)\mathrm{T}} \tag{32}$$

Here, $\mathrm{diag}(\boldsymbol{A})$ denotes a diagonal matrix composed of the diagonal entries of the matrix $\boldsymbol{A}$. The construction of $\mathrm{diag}\left(\boldsymbol{S}^{(k)}\right)$ is, however, not straightforward because the local Schur complement $\boldsymbol{S}^{(k)}$ is not constructed explicitly. Therefore, $\boldsymbol{S}^{(k)}$ is approximated by $\boldsymbol{K}_{\mathrm{BB}}^{(k)}$ as follows:

$$\boldsymbol{M}_{\mathrm{Diag\text{-}KBB}} = \sum_{k=1}^{N_{\mathrm{K}}} \boldsymbol{N}^{(k)} \mathrm{diag}\left( \boldsymbol{K}^{(k)} \right) \boldsymbol{N}^{(k)\mathrm{T}} + \sum_{k=1}^{N_{\mathrm{B}}} \boldsymbol{N}_{\mathrm{B}}^{(k)} \mathrm{diag}\left( \boldsymbol{K}_{\mathrm{BB}}^{(k)} \right) \boldsymbol{N}_{\mathrm{B}}^{(k)\mathrm{T}} \tag{33}$$

Even for the subdomains replaced by ROMs, $\boldsymbol{K}_{\mathrm{BB}}^{(k)}$ for the HDM must be constructed because it is used in Eq. (27). Therefore, $\boldsymbol{M}_{\mathrm{Diag\text{-}KBB}}$ in Eq. (33) can be computed easily. Accordingly, $\boldsymbol{M}_{\mathrm{Diag\text{-}KBB}}^{-1}$ is employed as the preconditioning matrix $\boldsymbol{M}^{-1}$ in Algorithm 1.

## 6. Generation of ROM Snapshots for the HDM–ROM Method

In the HDM–ROM method, the interfaces between a ROM subdomain and the other subdomains are treated as Dirichlet boundaries. In the present study, the ROM subdomain is assumed to have no boundaries other than the interface and the fixed boundary (i.e., a Dirichlet boundary with zero prescribed displacement); that is, neither Neumann boundaries for applied loads nor Dirichlet boundaries with prescribed nonzero displacements are considered. In addition, no body forces are applied. Under these assumptions, snapshot generation methods for constructing the ROB $\boldsymbol{V}$ are investigated.

The snapshots are computed by applying prescribed displacements to the nodes on the interface between subdomains. Three methods for applying the prescribed displacements are considered: (1) A unit displacement (unity is sufficient because the problem is linear) is applied to each degree of freedom (DOF) on the interface, while zero displacement is applied to all the other DOFs. (2) When the interface is planar, the prescribed displacements at the interface nodes are computed from the displacement fields obtained by tilting the plane in two directions while maintaining its planarity. (3) The entire analysis domain is modeled using an HDM and analyzed with the parallel CG solver, and the intermediate solutions obtained during the CG iterations are collected as snapshots.

In Method (1), the solution corresponding to an arbitrary displacement field prescribed on the interface can be represented as a linear combination of the snapshots. Consequently, the ROM constructed by this method can reproduce the HDM exactly if the number of ROB vectors is not reduced by singular value decomposition. Method (2) is expected to be effective, for example, when a component in an assembly structure is replaced by a ROM and connected to the other components through a supporting structure with a small cross-sectional area. In this case, the interface between the supporting structure and the component becomes the interface between subdomains, and the interface may be assumed to remain planar after deformation. In Method (3), reducing the number of ROB vectors below the number of collected snapshots by singular value decomposition is considered meaningful.

## 7. Implementation

The method proposed in the present paper was implemented using the ADVENTURE system (ADVENTURE Project, 2026), which consists of a set of open-source software modules developed for large-scale parallel computational mechanics. However, the parallel CG solver used as the basis for the implementation of the HDM–ROM method is not publicly available and is different from ADVENTURE_Solid, which is a parallel finite element analysis code based on hierarchical domain decomposition and is one of the modules of the ADVENTURE system. Nevertheless, the parallel CG solver employs libfem, which is the finite element library of ADVENTURE_Solid, for finite element computations. In addition, ADVENTURE_IO, which is also a module of the ADVENTURE system, is employed for input and output operations. Furthermore, ADVENTURE_Metis, a domain decomposition module based on the graph partitioning codes Metis and ParMetis, can decompose input data into hierarchical or single-layered domain decomposed data. Note that a non-overlapping domain decomposition is adopted. The parallel CG solver uses the single-layered decomposed data. This partitioning is performed using only ParMetis. The solver is parallelized using MPI, and one subdomain is assigned to one MPI process.

In the proposed method, the region to be modeled by a ROM must be treated as a single subdomain. In the case of multiple ROMs, each ROM is treated as a single subdomain. The remaining regions to be modeled by HDM must be automatically divided into an appropriate number of subdomains, as in a conventional domain decomposition method. Therefore, the HDM region is first divided into subdomains using ADVENTURE_Metis, and subdomain IDs are assigned on the basis of the resulting partitioning. Different subdomain IDs are assigned to the ROM regions. The entire mesh with the embedded subdomain IDs is then redivided using a modified version of ADVENTURE_Metis, in which domain decomposition is performed on the basis of the embedded subdomain IDs instead of ParMetis.

The construction of ROMs requires singular value decomposition of dense matrices and the use of ROMs requires the solution of linear systems with dense coefficient matrices. The computation in Eq. (30) to recover the HDM solution also requires dense matrix operations. These computations are implemented on GPUs using cuSOLVER and cuBLAS (NVIDIA, 2026), which are called only in the MPI processes that handle the ROMs in the parallel CG solver. They are also implemented on CPUs using LAPACKE in Intel oneAPI MKL (Intel, 2026).

## 8. Numerical Examples

The first analysis model is a rectangular solid discretized by linear hexahedral elements. The numbers of elements in the three coordinate directions are $3 \times 3 \times 12$, resulting in 208 nodes and 108 elements. The domain decomposed mesh and boundary conditions are shown in Fig. 3. The bottom surface is fixed, and a horizontal prescribed displacement is applied to the top surface. The bottom cubic region is modeled by a ROM and is treated as a single subdomain. The remaining HDM region is divided into six subdomains. The domain decomposition is performed using the method described in Section 7.

First, the three snapshot generation methods described in Section 6 are examined. Since the bottom surface of the ROM region is fixed in this analysis model, it is also fixed during the snapshot generation. The top surface of the ROM region forms the interface between subdomains, and snapshots are computed by applying prescribed displacements according to the three methods. The number of nodes on the interface of the ROM region is $4 \times 4 = 16$. In Method (1), prescribed displacements are applied in the three coordinate directions at each node, and therefore $16 \times 3 = 48$ snapshots are generated. Since the deformation modes represented by these snapshots are obviously linearly independent, the number of ROB vectors is not reduced on the basis of the singular values. The ROM constructed by Method (1) is referred to as **ROM_point**. Examples of the snapshots are shown in Fig. 4.

In Method (2), the interface plane is tilted in four patterns as shown in Fig. 5, and four snapshots are calculated by the prescribed displacements determined from the four tilted planes. Also in this case, the number of ROB vectors is not reduced on the basis of the singular values. The ROM constructed by Method (2) is referred to as **ROM_plane**.

In Method (3), the entire analysis model is first analyzed as an HDM, and the nodal displacement vectors obtained in the ROM region during the CG iterations are collected as snapshots. The ROM constructed by Method (3) is referred to as **ROM_cg.** Since the CG method converges in 53 iterations, 53 snapshots are collected.

The number of ROB vectors is determined using Eq. (7) with $\eta = 0.9999$ and $0.99999$. The numbers of ROB vectors are 20 and 24, respectively. The resulting ROMs are referred to as **ROM_cg_V20** and **ROM_cg_V24**, respectively. Examples of the snapshots are shown in Fig. 6.

Hereafter, an HDM–ROM model is referred to by the name of the ROM used in the model; for example, an HDM–ROM model using ROM_point is referred to as the HDM–ROM_point model. The maximum relative errors in the infinity norm for the HDM–ROM models with respect to the displacement vector obtained by the full HDM model are calculated for each subdomain to investigate the distribution of the errors and listed in Table 1.

The HDM–ROM_point model reproduced exactly the same solution as the HDM model. The HDM–ROM_plane model approximately reproduced the HDM solution. The HDM–ROM_cg_V24 model reproduced results that were almost identical to those of the HDM model. For the HDM–ROM_cg_V20 model, however, the maximum relative error was 2.24% in the HDM region and 10.30% in the ROM region. Note also that the HDM solution could not be reproduced for analysis under the $y$-direction load when the HDM–ROM_cg model constructed from snapshots generated only for the $x$-direction load was employed. On the other hand, the HDM–ROM_point model reproduced exactly the same solution as the HDM model for various load conditions, except when the load was applied to the ROM region. The HDM–ROM_plane model also approximately reproduced the HDM solution for various load conditions.

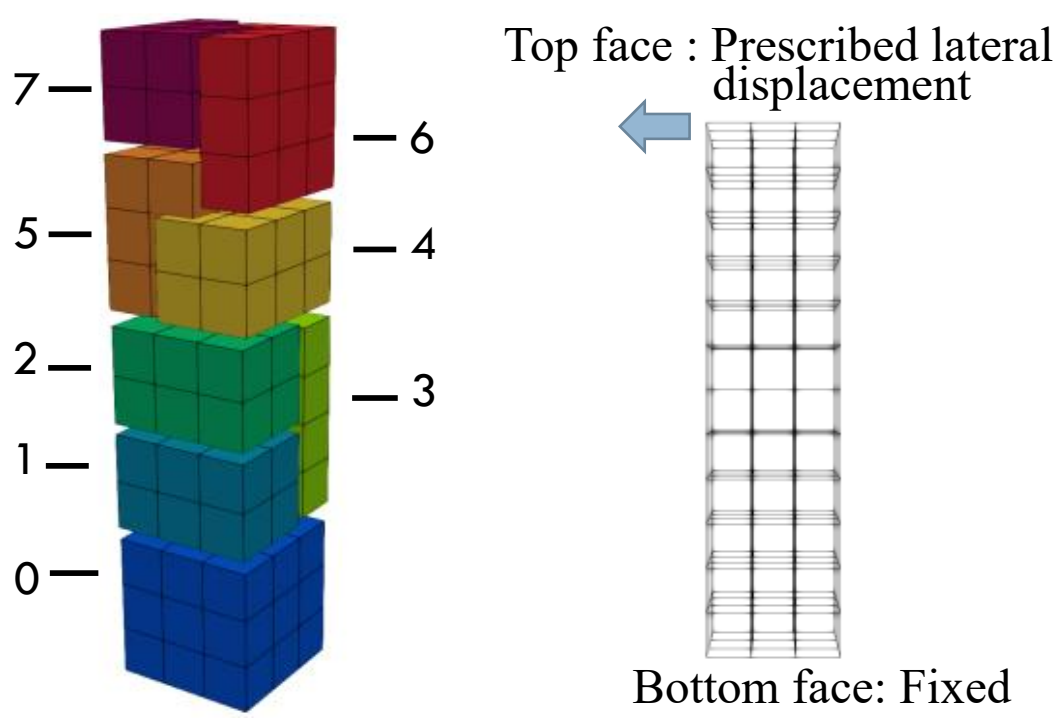

Fig. 3 HDM–ROM model for rectangular solid (domain decomposition and boudary conditions).

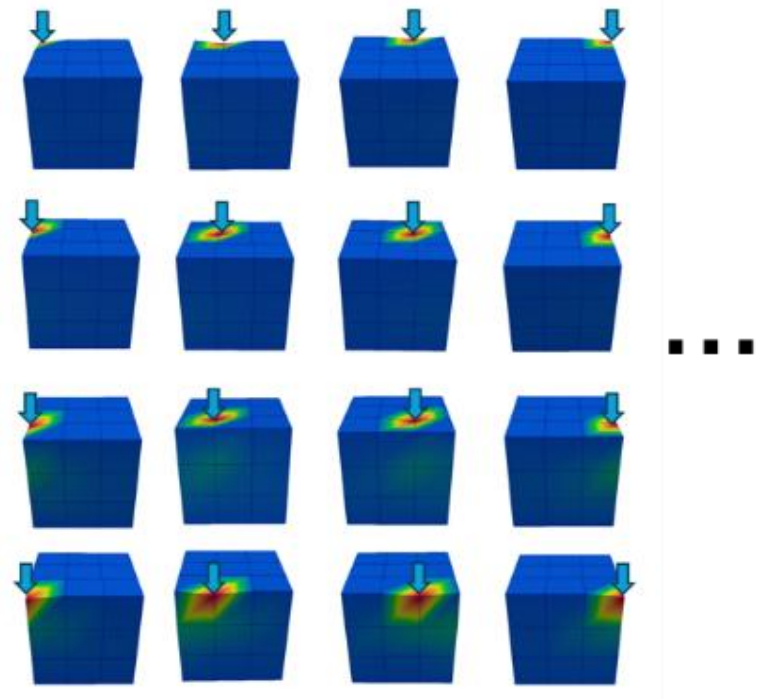

Fig. 4 Examples of snapshots for ROM_point model.

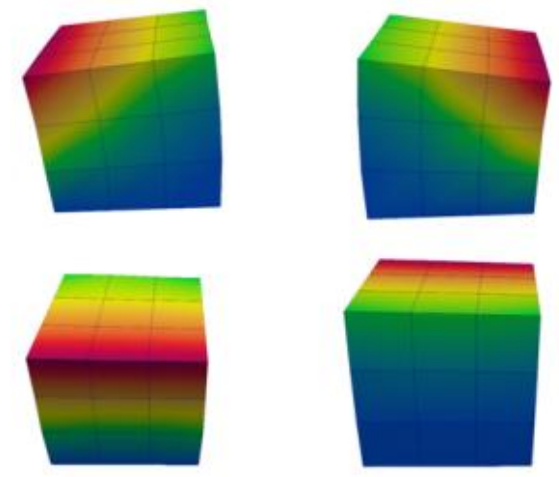

Fig. 5 Snapshots for ROM_plane model.

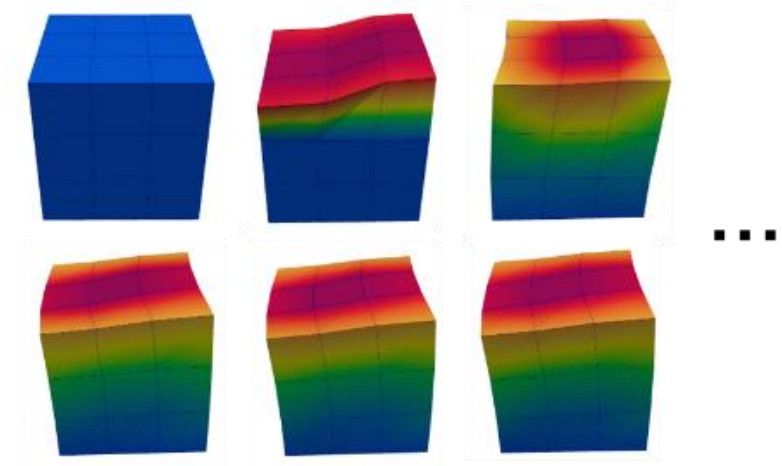

Fig. 6 Examples of snapshots for ROM_CG model.

Table 1 The maximum relative errors for the HDM–ROM models with respect to the displacement vector obtained by the full HDM model for each subdomain.

| Subdomain | HDM–ROM_point | HDM–ROM_plane | HDM–ROM_cg_V20 | HDM–ROM_cg_V24 |
|---|---|---|---|---|
| 0 | 0.00% | 9.14% | 10.30% | 0.00% |
| 1 | 0.00% | 3.52% | 2.24% | 0.00% |
| 2 | 0.00% | 1.62% | 1.24% | 0.00% |
| 3 | 0.00% | 1.94% | 1.24% | 0.00% |
| 4 | 0.00% | 0.80% | 0.72% | 0.00% |
| 5 | 0.00% | 0.68% | 0.61% | 0.00% |
| 6 | 0.00% | 0.32% | 0.31% | 0.00% |
| 7 | 0.00% | 0.21% | 0.22% | 0.00% |

The second analysis model is a simplified nuclear power plant model consisting of a building and four SMRs, as shown in Fig. 1. Boundary conditions are shown in Fig. 7. For the evaluation of computational performance of the GPU implementation, a model with one SMR is also constructed. The models with one SMR and four SMRs are referred to as SMR1 and SMR4, respectively. The geometry of the building in SMR1 is the same as that in SMR4. The geometry is defined using 3D CAD software so that each SMR and the building constituted a separate volume, and a linear tetrahedral mesh is generated from this geometry. The SMRs are modeled by ROMs. SMR4 is divided into eight subdomains, as shown in Fig. 8, whereas SMR1 is divided into five subdomains. Therefore, the number of subdomains in the HDM region is four in both models. Since the SMRs are modeled by ROMs, the domain decomposition is performed so that each SMR corresponds to a single subdomain.

Methods (1) and (2) described in Section 6 are used to construct ROMs for the SMRs. Therefore, four ROMs are constructed for SMR4. Note that the meshes of the SMRs are slightly different although the geometries defined by the 3D CAD software are identical. Therefore, four ROMs are constructed separately. As shown in Fig. 9, each SMR is supported on three surfaces. The three supporting surfaces of each SMR serve as the interfaces between subdomains. Table 2 lists the mesh information for the analysis models and for each SMR. In Method (1), prescribed displacements are applied in the three coordinate directions at each node. Therefore, the number of snapshots for SMR#0, for example, is $981 \times 3 = 2{,}943$. The number of ROB vectors is determined using Eq. (7) with $\eta = 0.8$ and $1.0$. The numbers of ROB vectors for SMR#0 are 1,501 and 2,943, respectively. The resulting ROMs are referred to as **ROM_point_80%** and **ROM_point_100%**, respectively. In Method (2), four Dirichlet boundary conditions are prescribed on each supporting surface, as in the rectangular solid model. Note that when a nonzero Dirichlet boundary condition is prescribed on one supporting surface, the other supporting surfaces are fixed. Since each SMR has three supporting surfaces, the total number of snapshots for one SMR is $4 \times 3 = 12$. Examples of the snapshots are shown in Fig. 10.

Table 3 lists the maximum relative error in the infinity norm for each subdomain with respect to the results obtained by the full HDM model. The HDM–ROM_point_100% model reproduced exactly the same solution as the HDM model. For the HDM–ROM_point_80% model, the maximum relative errors in the displacement were 0.25 % in the HDM region and 0.63 % in the ROM regions. For the HDM–ROM_plane model, on the other hand, the maximum relative errors were 66.69% in the HDM region and 55.28% in the ROM regions. The convergence histories of the relative norm of the CG residual vector, normalized by its initial value, for SMR4 are shown in Fig. 10. Since the number of DOFs in the ROM regions to be solved by the CG method is smaller than that in the full HDM model, the HDM–ROM model converged faster than the full HDM model. The diagonal scaling preconditioner using Eq. (33) for the HDM–ROM model also improved the convergence.

The computation times and the numbers of CG iterations required for convergence are summarized in Table 4. The diagonal scaling preconditioner was used in all cases. The computations were performed on a workstation equipped with an AMD EPYC 9554 CPU (64 cores / 128 threads) and an NVIDIA A100 GPU (432 Tensor Cores and 3456 FP64 CUDA cores). One subdomain was assigned to one MPI process, and each MPI process was assigned to one CPU core. The domain decomposition for the HDM is the same as that for the corresponding HDM–ROM model; that is, each SMR is treated as a single subdomain. Therefore, this decomposition is not optimal for the HDM. The ROM computations were performed using LAPACKE in Intel MKL on CPUs for SMR4, whereas those for SMR1 were performed using either

LAPACKE or cuSOLVER and cuBLAS on a GPU. One thread was used for the LAPACKE implementation in each MPI process. Note that the computation times for data reading, preparation of arrays used for the ROM, and factorization of the ROM coefficient matrix in Eq. (9) are excluded from the total computation time. It is also noted that the factorization time is very short, approximately 242 ms with LAPACKE and 5 ms with cuSOLVER for ROM_point_100% of one SMR. The computation times per CG iteration for SMR1 and SMR4 using five and eight MPI processes, respectively, were almost identical, indicating that the computations for the four ROMs in SMR4 were successfully performed in parallel. In addition, the computation time per CG iteration was longer than that for the full HDM model because the dense matrix operations required for the ROMs are computationally expensive. Nevertheless, the total computation time for the HDM–ROM model was shorter than that for the full HDM model because the number of CG iterations was reduced. It is noteworthy that the computation time for the ROM computations was drastically reduced by the GPU implementation. A more detailed evaluation of the computational performance, together with the details of the hybrid MPI-CUDA implementation, will be presented in a future paper.

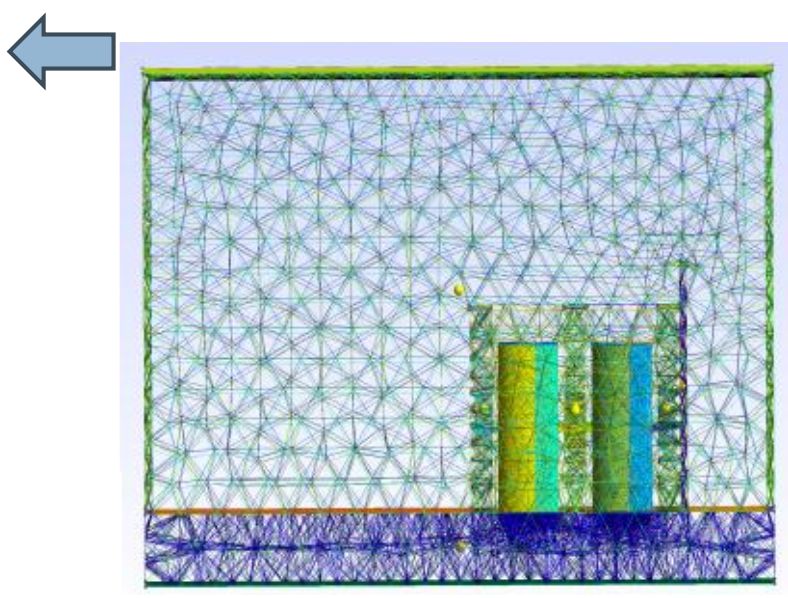

Fig. 7 Boundary conditions for the nuclear power plant model.

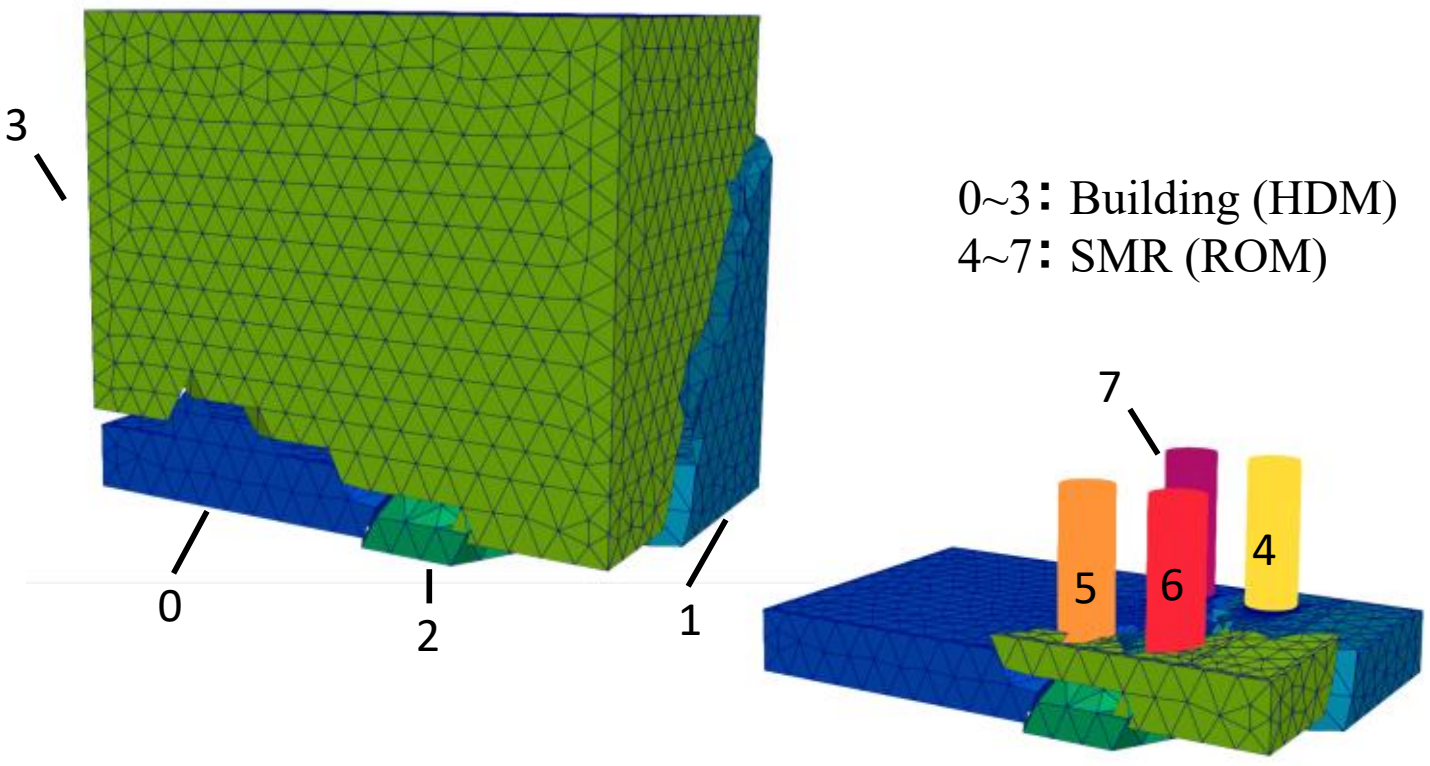

Fig. 8 Domain decomposition for the nuclear power plant model. Numbers denotes subdomain number.

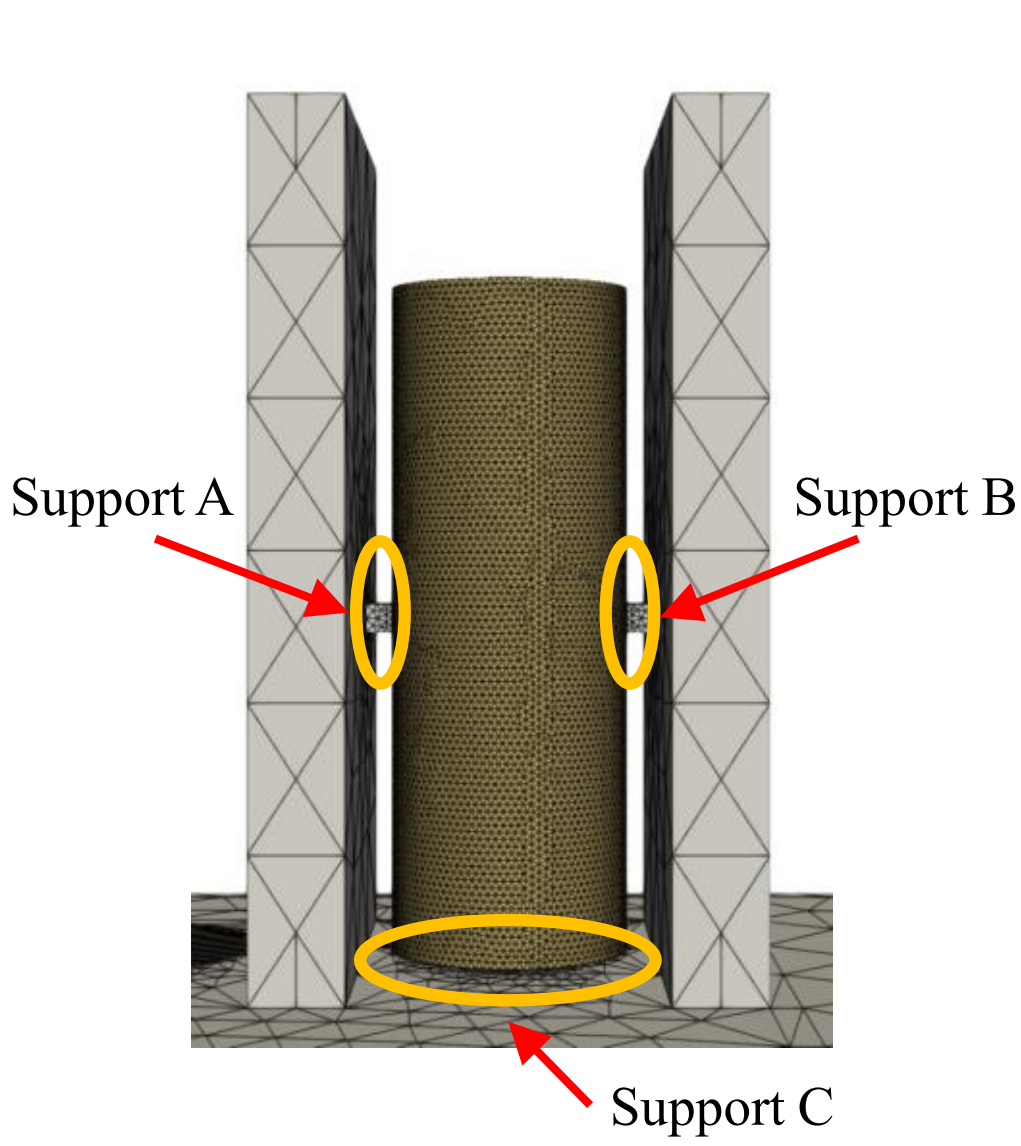

Fig. 9 Magnified view of an SMR mesh with support structures circled with yellow ellipses.

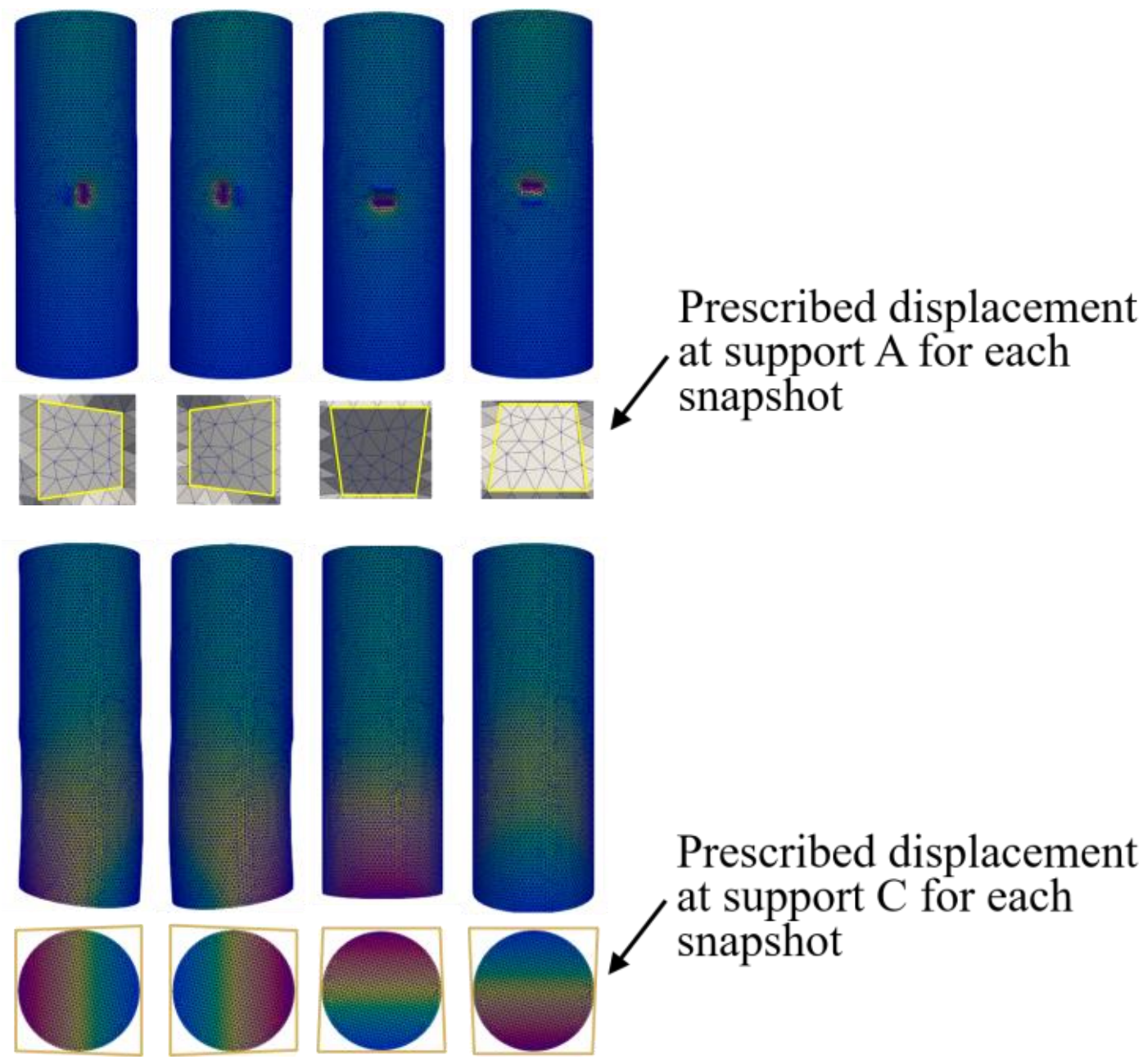

Fig. 10 Examples of snapshots for a ROM of an SMR (the upper row shows snapshots with respect to the upper support structure, and the lower row shows those with respect to the lower support structure).

Table 2 Mesh specifications.

| Model | Entire model /Each SMR | Number of elements | Number of nodes | Total number of nodes on supporting surfaces |
|---|---|---|---|---|
| SMR4 | Entire | 1,079,454 | 203,769 | – |
| | SMR #0 | 261,482 | 49,388 | 981 |
| | SMR #1 | 261,869 | 49,441 | 981 |
| | SMR #2 | 261,688 | 49,437 | 980 |
| | SMR #3 | 261,960 | 49,456 | 981 |
| SMR1 | Entire | 286,902 | 56,363 | – |
| | SMR | 259,697 | 49,092 | 980 |

Table 3 The maximum relative errors for the HDM–ROM models with respect to the displacement vector obtained by the full HDM model for each subdomain.

| Subdomain | HDM–ROM_point_100% | HDM–ROM_point_80% | HDM–ROM_plane |
|---|---|---|---|
| 0 | 0.00% | 0.07% | 4.09% |
| 1 | 0.00% | 0.06% | 2.64% |
| 2 | 0.00% | 0.25% | 66.69% |
| 3 | 0.00% | 0.02% | 1.67% |
| 4 | 0.00% | 0.63% | 50.30% |
| 5 | 0.00% | 0.30% | 53.47% |
| 6 | 0.00% | 0.39% | 50.18% |
| 7 | 0.00% | 0.63% | 55.28% |

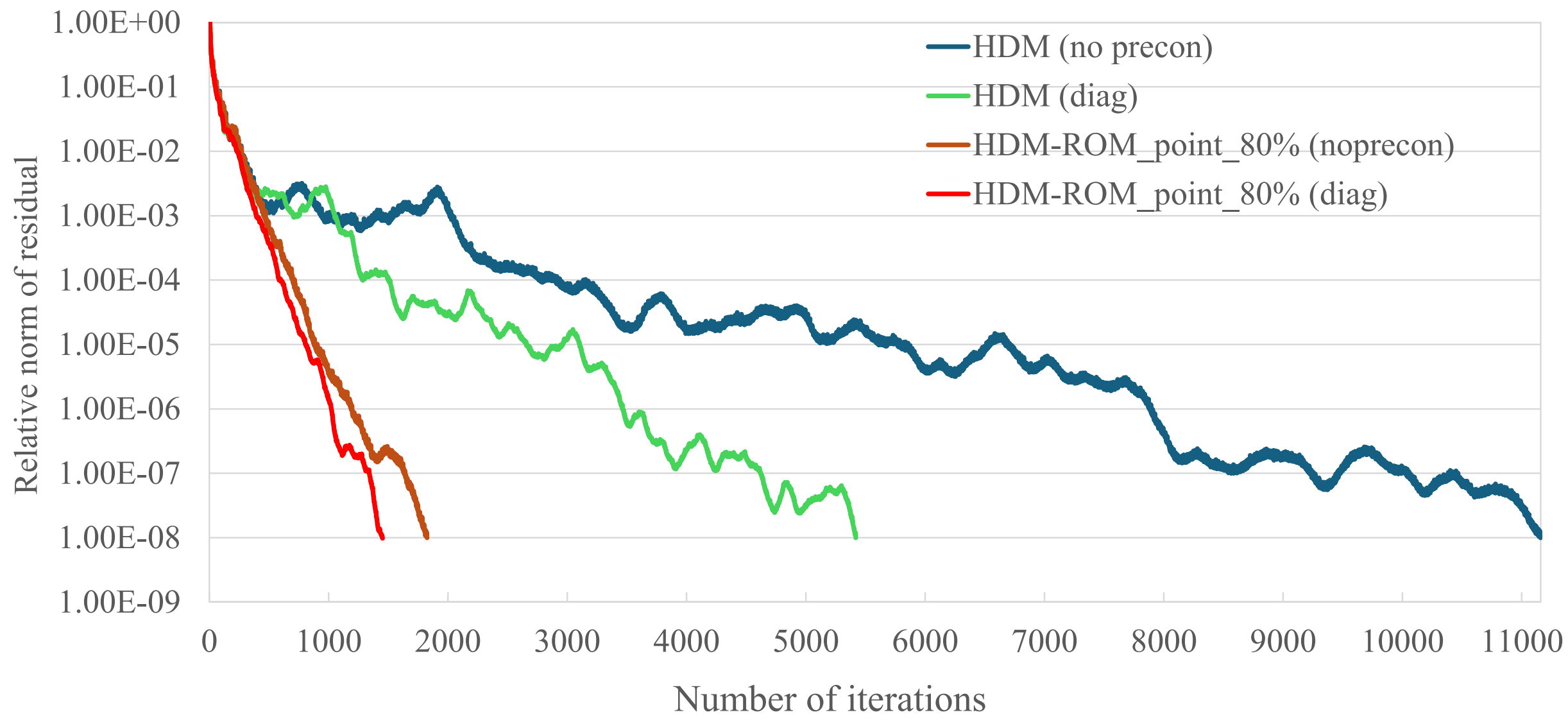

Fig. 11 Convergence histories of the CG method (SMR4).

Table 4 Computation times. The diagonal scaling preconditioner was used in all cases.

| Case | CPU/GPU | Number of iterations | Computation time (s) | Computation time for one CG iteration (s) |
|---|---|---|---|---|
| SMR4 HDM | CPU | 5426 | 76.7 | 0.0142 |
| SMR1 HDM | CPU | 3481 | 48.5 | 0.0139 |
| SMR4 HDM–ROM_point_100% | CPU | 1452 | 371.3 | 0.2558 |
| SMR4 HDM–ROM_point_80% | CPU | 1451 | 191.3 | 0.1318 |
| SMR4 HDM–ROM_plane | CPU | 1435 | 4.6 | 0.0032 |
| SMR1 HDM–ROM_point_100% | CPU | 1944 | 404.1 | 0.2079 |
| SMR1 HDM–ROM_point_100% | GPU | 1947 | 13.0 | 0.0067 |
| SMR1 HDM–ROM_point_80% | CPU | 1942 | 206.3 | 0.1063 |
| SMR1 HDM–ROM_point_80% | GPU | 1943 | 8.4 | 0.0043 |
| SMR1 HDM–ROM_plane | CPU | 1866 | 4.4 | 0.0023 |
| SMR1 HDM–ROM_plane | GPU | 1866 | 3.2 | 0.0017 |

## 9. Conclusions

In the present paper, the HDM–ROM method is proposed; that is, HDM and multiple projection-based ROMs are combined within the framework of an iterative DDM. The CG method is employed as the iterative solver. The HDM domain is subdivided into multiple subdomains and treated using a non-substructuring DDM, whereas the multiple ROM domains are treated using a substructuring-based DDM, in which the static condensation procedure is approximated using ROMs. The proposed method is implemented within the framework of an existing DDM-based iterative solver. The ROM computations are implemented on GPUs using CUDA. The following findings are obtained.

- When the ROB was constructed from snapshots computed by applying prescribed unit displacements to each DOF on the interface while fixing all the other DOFs (ROM_point model), the HDM–ROM method reproduced exactly the same solution as the HDM model.
- When the ROB was constructed from snapshots computed by tilting the interface plane in two directions while maintaining its planarity (ROM_plane model), the HDM–ROM model approximately reproduced the HDM solution only when the deformed interface could be approximated by a planar surface.
- When the ROB was constructed from the nodal displacement vectors obtained during the CG iterations for the entire analysis model, the ROM using the ROB could not reproduce the HDM solution under a load condition different from that used for snapshot collection.
- The CG method exhibited good convergence for the HDM–ROM model. The approximate diagonal scaling preconditioner reduced the number of iterations required for convergence.
- The HDM–ROM method was applied to a power plant model with SMRs, which is an example of an assembly structure with a complex geometry. The solver converged successfully, and the HDM–ROM_point model reproduced the exact solution. However, the HDM–ROM_plane model could not reproduce the solution with sufficient accuracy, although the solver also converged successfully.
- The GPU implementation significantly reduced the computation time..

**Acknowledgements**

This study was supported by JSPS KAKENHI (Grant Number 24K14984) and partly supported by the Kajima Foundation's General Research Grants. This work used the computational resources of the supercomputer Fugaku provided by the RIKEN Center for Computational Science (Project ID: hp250298) for the snapshot computations. Valuable advice was provided by Dr. Takuzo Yamashita (NIED). Sincere appreciation is expressed for these support and advice.